%% file: gasper_trebels4.tex
\documentstyle[12pt,dina4]{article}
 \input{macrola2.tex}

\begin{document}
 \begin{center}
 {\LARGE{}Ultraspherical multipliers revisited }\\[.4cm]
 {\sc{}George  Gasper\footnote{Department of Mathematics, Northwestern 
 University, Evanston, IL 60208, USA. The work of this author was 
 supported in part by the National Science   Foundation under grant 
 DMS--9401452.}
 and Walter  Trebels\footnote{Fachbereich Mathematik, TH Darmstadt, 
 Schlo\ss{}gartenstr.7, D--64289
 Darmstadt, Germany.} }\\[.4cm]
 {\it Dedicated to L\'aszl\'o Leindler on the occasion of his 60th birthday} \\ 
 [.4cm]
 {(March. 28, 1995 version)}
 \end{center}

 \bigskip
 {\bf Abstract.}  Sufficient ultraspherical multiplier criteria are refined 
 in  such a way  that they are comparable with necessary multiplier conditions.
 Also new necessary conditions for Jacobi multipliers are 
 deduced which, in particular, imply known Cohen type inequalities. 
Muckenhoupt's 
 transplantation theorem is used in an essential way.


\bigskip
{\bf AMS(MOS) subject classifications.} 33C45, 42A45, 42C10 

\bigskip 
\section{Introduction}
Quite sharp sufficient conditions for ultraspherical multipliers are contained 
in papers by Muckenhoupt and Stein \cite{must}, Bonami and Clerc \cite{bocl}, 
Connett and Schwartz \cite{cosc}, Gasper and Trebels \cite{jamu} and 
Muckenhoupt \cite{mu}. In \cite{hausdorff} 
we gave comparable necessary conditions for Jacobi multipliers with
parameters $(\alpha ,-1/2)$ in the ``natural'' weight case 
(see \cite[p. 2]{mu} and below). It is the goal of this paper 
to develop necessary conditions for ultraspherical (Jacobi) multipliers and to 
weaken the sufficient ones in such a way that they are comparable with the 
necessary ones.    This is done by decomposing the relevant functions 
into even and odd parts; thus, by the quadratic transformations in 
\cite[(4.1.5)]{szego}, reducing the problem of controlling the multiplier 
sequence $\{ m_k\} $ to a discussion of the 
subsequence $\{ m_{2k} \} $ in the Jacobi case $(\alpha , -1/2)$ with natural 
weight and of $\{ m_{2k+1} \} $ for the parameters $(\alpha ,1/2)$ (with an 
additional weight). An essential tool is Muckenhoupt's \cite{mu} 
transplantation theorem.

\medskip \noindent
To become more precise let us introduce some notation. In view of the above it 
is reasonable to work within the framework of Jacobi expansions --- the 
conversion to the standard notation for ultraspherical polynomials given in 
 \cite[(4.7.1)]{szego} does not change the involved multiplier 
spaces.\\
Fix $\alpha \ge \beta \ge -1/2$ and let $L^p_{(a ,b)}\, , 1\le p < \infty \, ,$ 
denote the space of measurable functions on $[0,\pi ]$ with finite norm 
$$\| f \| _{L^p_{(a ,b)}} =\Big( \int _0^\pi \Big| f(\theta )\Big| ^p 
\Big( \sin \frac{\theta }{2} \Big) ^{2a+1} 
\Big( \cos \frac{\theta }{2} \Big) ^{2b+1}d\theta  \Big) ^{1/p} .$$
If $a=b$ we use the abbreviation $L^p_a = L^p_{(a,a)}$.
The ``natural'' weight case for expansions in Jacobi polynomials (when 
there is a nice convolution structure) is the case
when $a=\alpha ,\; b=\beta $. Define the normalized Jacobi polynomials by 
$R_k^{(\alpha ,\beta )} (x)=P_k^{(\alpha ,\beta )} (x)/P_k^{(\alpha ,
\beta )} (1)$, where $P_k^{(\alpha ,\beta )} (x)$ is the Jacobi polynomial of 
degree $k$ and order $(\alpha ,\beta )$, see \cite{szego}. For  $f\in 
L^1_{(\alpha ,\beta )} $, its $k$-th
Fourier--Jacobi coefficient ${\hat f}_{(\alpha ,\beta)} (k)$ is defined by
$${\hat f}_{(\alpha ,\beta)} (k)=\int _0^\pi f(\theta )R_k^{(\alpha ,\beta )} 
(\cos \theta )\Big( \sin \frac{\theta }{2} \Big) ^{2\alpha +1} 
\Big( \cos \frac{\theta }{2} \Big) ^{2\beta +1} d\theta \, .$$
Then $f$ has an expansion of the form
$$f(\theta ) \sim \sum _{k=0}^\infty {\hat f}_{(\alpha ,\beta)} (k)
h_k^{(\alpha ,\beta )} R_k^{(\alpha ,\beta )} (\cos \theta )\, ,$$
where the normalizing factors $h_k^{(\alpha ,\beta )} $ are given by 
$h_k^{(\alpha ,\beta )} = 
\| R_k^{(\alpha ,\beta )} (\cos \theta )\|_{L^2_{(\alpha ,\beta )}}
 ^{-2}\approx (k+1)^{2\alpha +1}$ (here the $\approx $ sign means that there are 
positive constants $C,C'$ such that $C'h_k^{(\alpha ,\beta )} \le (k+1)
^{2\alpha +1} \le C h_k^{(\alpha ,\beta )}$ holds).

\medskip \noindent
A sequence $m=\{ m_k\} _{k=0}^\infty \in l^\infty $ is called a multiplier on 
$L^p_{(a,b)}$ with respect to an expansion into Jacobi polynomials of order 
$(\alpha ,\beta )$, notation $m\in M^p_{(\alpha ,\beta );(a,b)}\, $, if for each 
$f\in L^p_{(a,b)}$ there exists a function $T_mf \in L^p_{(a,b)}$ with
\begin{equation}\label{defmult}
 T_mf(\theta ) \sim \sum _{k=0}^\infty m_k {\hat f}_{(\alpha ,\beta)} (k)
h_k^{(\alpha ,\beta )} R_k^{(\alpha ,\beta )} (\cos \theta )\, ,\quad 
\| T_mf \| _{L^p_{(a,b)}} \le C \| f\| _{L^p_{(a,b)}}\, .
\end{equation}
The smallest constant $C$ independent of $f$ for which this holds is called the 
multiplier norm of $m$ and is denoted by $\| m\| _{M^p_{(\alpha ,\beta 
);(a,b)}}$.  If $\alpha = \beta$ and $a=b,$ we write $M^p_{\alpha ;a}.$

\bigskip \noindent
Now decompose a function $f\in L^p_\alpha $ into its even part 
$f_e$ and its odd part $f_o$ with respect to the line $\theta =\pi /2$:
$$ f_e(\theta)=\{ f(\theta)+ f(\pi -\theta) \}/2 \, , 
\quad f_o=f-f_e\, .$$
Obviously, this decomposition is unique and 
there holds for the Fourier--Jacobi coefficients (observe that $R_k^{(\alpha 
,\alpha )}(x)$ is even when $k$ is even and odd when $k$ is odd)
\begin{equation}\label{decompcoeff}
 (f_e)\,  {\hat {} }_{(\alpha ,\alpha )} (k)=
     \left\{ \begin{array}{l@{\; ,\; }l}         
{\hat f}_{(\alpha ,\alpha )} (k)  & k \; {\rm even} \\
                     0 & k \; {\rm odd}  \end{array} \right.  ,\quad 
(f_o)\,  {\hat {} }_{(\alpha ,\alpha )} (k)=
     \left\{ \begin{array}{l@{\; ,\; }l}             
0 & k \; {\rm even} \\
                {\hat f}_{(\alpha ,\alpha )} (k)  & k \; {\rm odd.}
\end{array}  \right. 
\end{equation}
Furthermore,
\begin{equation}\label{decompnorm}
\| f \| _{L^p_\alpha } \le \| f_e \| _{L^p_\alpha }
+\| f_o \| _{L^p_\alpha } \le 2 \| f \| _{L^p_\alpha } \, ,\; \; 
1\le p < \infty \, .
\end{equation}
In particular, the uniqueness theorem shows that 
\begin{equation}\label{decompsereven}
f_e(\theta ) \sim \sum _{k=0}^\infty {\hat f}_{(\alpha ,\alpha )} (2k)
h_{2k}^{(\alpha ,\alpha )} R_{2k}^{(\alpha ,\alpha )} (\cos \theta )\, ,
\end{equation}
\begin{equation}\label{decompserodd}
f_o(\theta ) \sim \sum _{k=0}^\infty {\hat f}_{(\alpha ,\alpha )} (2k+1)
h_{2k+1}^{(\alpha ,\alpha )} R_{2k+1}^{(\alpha ,\alpha )} (\cos \theta )\, .
\end{equation}

\medskip \noindent
Given a sequence $m=\{ m_k\} $ it is clear by the above that its $M^p_{\alpha 
;\alpha }$--multiplier norm is equivalent to the 
multiplier norm of $ m$ restricted to the subspace of even $L^p_\alpha 
$--functions (with respect to the line $\theta =\pi /2$) plus  
the multiplier norm of $ m$ restricted to the subspace of odd $L^p_\alpha 
$--functions, i.e., 
\begin{equation}\label{decompmult}
 \| m \|_{M^p_{\alpha ;\alpha }}
\approx \| m \|_{M^p_{\alpha ;\alpha } \big| _{{\rm even}}}
+ \| m \|_{M^p_{\alpha ;\alpha } \big| _{{\rm odd}}}\, .
\end{equation}
We can now state our first theorem.

\medskip \noindent
\thm{
Assume $\alpha \ge -1/2$ and define subsequences $m_e$ and $m_o$ of a given 
sequence $m$ by $(m_e)_k =m_{2k},\; (m_o)_k=m_{2k+1},\; k\in {\bf N}_0$.
\begin{itemize}
\item[a)] If $1 \le p < \infty, $ then there holds
$$\|  m \|_{M^p_{\alpha ;\alpha } \big| _{{\rm even}}}
\approx \| m_e\| _{M^p_{(\alpha ,-1/2);(\alpha ,-1/2)}}\, ,$$
$$\|  m \|_{M^p_{\alpha ;\alpha } \big| _{{\rm odd}}}
\approx \| m_o \| _{M^p_{(\alpha ,1/2);(\alpha ,(p-1)/2)}}\, ,$$
whenever one side in each of the equivalences is finite.
\item[b)] If $1<p<\infty, $ then
$$ \| m \|_{M^p_{\alpha ;\alpha }} \approx 
\| m_e\| _{M^p_{(\alpha ,-1/2);(\alpha ,-1/2)}} +
\| m_o\| _{M^p_{(\alpha ,-1/2);(\alpha ,-1/2)}} \, .$$
\end{itemize}
} 

\medskip \noindent
We will combine Theorem 1.1 with known sufficient criteria and 
necessary ones. To this end  
define the fractional difference operator of order $\mu ,\; \mu \in {\bf R},$ 
with increment $\kappa \in {\bf N}$ by 
$$\Delta _\kappa ^\mu m_k=\sum _{j=0}^\infty A_j^{-\mu -1} m_{k+\kappa j}, \quad
A_j^\mu =\frac{ \Gamma (j+\mu +1)}{\Gamma (j+1) \Gamma (\mu +1)} \, ,$$
whenever the series converges; when $\kappa =1$ we write $\Delta ^\mu  =\Delta 
_1^\mu .$ 
An application of the multiplier criteria from \cite[Theorem 4]{jamu},  
\cite[(3.8)]{hausdorff} as well as of Askey's \cite{as} Marcinkiewicz multiplier 
theorem for Jacobi expansions 
to the sequences $ m_e$ and $ m_o$ and the observation that 
$$\| m\| ^q_\infty  + \sup _{N\in {\bf N}_0} 
 \sum _{k=N}^{2N}\big| (k+1)^\mu \Delta _2^\mu m_k \big| ^q \frac{1}{k+1} 
\approx \quad \quad \quad \quad \quad \quad $$
$$ \| m\| ^q_\infty  + \sup _{N\in {\bf N}_0}  \sum _{k=N}^{2N}\big| (k+1)^\mu  
\Delta ^\mu  (m_e)_k \big| ^q \frac{1}{k+1} +
\sup _{N\in {\bf N}_0} \sum _{k=N}^{2N}\big| (k+1)^\mu  
\Delta ^\mu (m_o)_k \big| ^q \frac{1}{k+1} $$
immediately lead to 

\medskip \noindent
\coro{
Let $\alpha \ge -1/2$, $1<p<\infty $, and let $\{ m_k \} \in l^\infty $ be a 
given sequence. 
\begin{itemize}
\item[a)] If $m$ satisfies for $\mu > {\rm max} \{ (2\alpha +2)|1/p-1/2|,1/2\}$ 
the condition
$$\| m\| _\infty  + \sup _{N\in {\bf N}_0} \Big( \sum _{k=N}^{2N}\big| 
(k+1)^\mu \Delta _2^\mu m_k \big| ^2 \frac{1}{k+1} \Big) ^{1/2} 
\le D < \infty \, ,$$
then $m\in M^p_{\alpha ;\alpha }$ and
$ \| m\| _{M^p_{\alpha ;\alpha }} \le C\, D.$
\item[b)] If $m\in M^p_{\alpha ;\alpha }\, ,$  then
$$ \| m\| _\infty +\sup_{N\in {\bf N}_0} \Big( \sum _{k=N}^{2N} \big| (k+1)^\nu  
\Delta _2^\nu m_k \big| ^{p'} \frac{1}{k+1} \Big) ^{1/p'} \le C
\| m\| _{M^p_{\alpha ;\alpha }}\, , $$
where $\nu \le (2\alpha +1)|1/p-1/2|$ and  $1/p+1/p'=1.$
\item[c)] If $m$ satisfies the condition
$$\| m\| _\infty + \sup _{N\in {\bf N}_0} \left( \sum _{k=N}^{2N} | 
\Delta _2 m_k |\right)  \le D^* < \infty \, ,\; $$
then $m\in M^p_{\alpha ;\alpha }$ if 
$1\le (4\alpha +4)/(2\alpha +3) <p<(4\alpha +4)/(2\alpha +1) \le \infty 
$ and $ \| m\| _{M^p_{\alpha ;\alpha }} \le C\, D^*.$
\end{itemize}
The constants $C$ in the above statements are independent of the 
sequences $m$.
}

\bigskip \noindent
{\bf Remarks} 1. The smoothness gap between the sufficient conditions in a) 
and the necessary ones in b) is essentially the gap which occurs in 
Sobolev embedding theorems (for the analogous result in the case of difference 
operators with increment $1$ see \cite[Theorem 5 b]{wbv}).

\medskip \noindent
2. We note that, in particular, Corollary 1.2 b) for half--integers 
$\alpha = (n-2)/2$ contains necessary 
conditions for zonal multipliers for spherical harmonic expansions (for the 
relevant notation and sufficient criteria see Strichartz \cite{stri}).

\medskip \noindent
3. Using the method which leads to Theorem 1.1 and its Corollary one can 
easily deduce weighted analogs. In the case $p=2$ and $c>0$ there holds
$$ \| m \|_{M^2_{\alpha ;\alpha +c}} \approx 
\| m_e\| _{M^2_{(\alpha ,-1/2);(\alpha +c,-1/2)}} +
\| m_o\| _{M^2_{(\alpha ,-1/2);(\alpha +c,-1/2)}} $$
provided the hypotheses of Muckenhoupt's transplantation theorem are satisfied, 
i.e., $(c+1-\alpha )/2$ is not a positive integer and the multipliers are 
defined on those subspaces of $L^p_{\alpha +c}$--functions for which 
${\hat f}_{(\alpha ,\alpha )} (k)=0,\; 0\le k \le {\rm max} \{ 0,\, 
[(c+1-\alpha )/2] \} -1.$ Now a characterization of multipliers for 
trigonometric series on weighted $L^2(-\pi ,\pi )$--spaces, due to Muckenhoupt, 
Wheeden and Young \cite[Theorems 10.1 and 10.2]{mwy}, 
(restricted to even functions) can be used to give
$$\| m\| _{M^2_{\alpha ;\alpha +c}} \approx 
\| m\| _\infty  + \sup _{N\in {\bf N}_0} \Big( \sum _{k=N}^{2N}\big| (k+1)^c 
\Delta _2^c m_k \big| ^2 \frac{1}{k+1} \Big) ^{1/2} ,$$
provided $c$ satisfies the condition $l+1/2 < c< l+3/2,\; l\in {\bf N}_0.$

\section{Proof of Theorem 1.1} 
\noindent
The reduction of ultraspherical multipliers to Jacobi multipliers with $\beta = 
-1/2$ or $1/2$ is accomplished by the transformation
formulas in \cite[(4.1.5)]{szego},
\begin{equation}\label{quadratic}
  R_{2k}^{(\alpha ,\alpha )} (\cos \theta )
=R_k^{(\alpha ,-1/2)} (\cos 2\theta ),\quad  
R_{2k+1}^{(\alpha ,\alpha )} (\cos \theta )
=\cos \theta \, R_k^{(\alpha ,1/2)} (\cos 2\theta ).
\end{equation}
The relevant Fourier--Jacobi coefficients are connected in the following way: 
\begin{equation}\label{transcoeffeven}
 2^{2\alpha +1}(f_e)\, {\hat {}}_{(\alpha ,\alpha)}(2k)
=[f_e(\theta /2)]\, {\hat {}}_{(\alpha ,-1/2)}(k)=:A_k,\quad 
k\in {\bf N}_0,
\end{equation}
\begin{equation}\label{transcoeffodd}
2^{2\alpha +1}(f_o)\, {\hat {}}_{(\alpha ,\alpha)}(2k+1)
=[ f_o(\theta /2)/\cos (\theta /2) ]\, 
{\hat {}}_{(\alpha ,1/2)}(k)=:B_k,\quad k\in {\bf N}_0.
\end{equation}
Furthermore, elementary computations give
\begin{equation}\label{normeven}
 h_{2k}^{(\alpha ,\alpha )} =2^{2\alpha +1} h_k^{(\alpha ,-1/2)},\quad
\| f_e\| _{L^p_\alpha } \approx 
\| f_e(\theta /2) \| _{L^p_{(\alpha ,-1/2)}} \, , \quad \quad \quad \quad 
\end{equation}
\begin{equation}\label{normodd}
 h_{2k+1}^{(\alpha ,\alpha )} =2^{2\alpha +1} h_k^{(\alpha ,1/2)} ,\quad 
\| f_o\| _{L^p_\alpha } \approx 
\| f_o(\theta /2)/(\cos (\theta /2))^{2/p} \| _{L^p_{(\alpha ,1/2)}} \, .
\end{equation}
This inserted in (\ref{decompsereven}) and (\ref{decompserodd}) leads for 
$f=f_e+f_o$ a cosine polynomial (i.e., a polynomial in powers of 
$\cos \theta $) to
\begin{equation}\label{transeven}
f_e(\theta )  = \sum _{k=0}^\infty A_k h_k^{(\alpha ,-1/2)} 
R_k^{(\alpha ,-1/2)}(\cos 2\theta ),
\end{equation}
\begin{equation}\label{transodd}
f_o(\theta ) = \sum _{k=0}^\infty B_k h_k^{(\alpha ,1/2)} 
\cos \theta \, R_k^{(\alpha ,1/2)}(\cos 2\theta ).
\end{equation}
Thus it follows from (\ref{transeven}) and (\ref{normeven}) that
\begin{eqnarray*}
\| T_m f_e\| _{L^p_\alpha } & \approx & \| \sum _{k=0}^\infty m_{2k} A_k 
h_k^{(\alpha ,-1/2)} \, R_k^{(\alpha ,-1/2)} (\cos 2\theta ) \| _{L^p_\alpha 
}\\
{} & \approx & \| \sum _{k=0}^\infty m_{2k} A_k 
h_k^{(\alpha ,-1/2)} \, R_k^{(\alpha ,-1/2)} (\cos \theta ) \| _{L^p_{(\alpha 
,-1/2)}}\\
{} & \le & C \|  m_e \| _{M^p_{(\alpha ,-1/2);(\alpha ,-1/2)}}
\| f_e (\theta /2) \| _{L^p_{(\alpha ,-1/2)}} \\
{} & \approx & 
\|  m_e \| _{M^p_{(\alpha ,-1/2);(\alpha ,-1/2)}}
\| f_e \| _{L^p_\alpha } 
\end{eqnarray*}
which implies
$$ \| \{ m_k\} \|_{M^p_{\alpha ;\alpha } \big| _{{\rm even}}} \le C 
\|  m_e \|_{M^p_{(\alpha ,-1/2);(\alpha ,-1/2)} } \, .$$
The converse is proved analogously by just starting with 
$\| \sum m_{2k} \dots \| _{L^p_{(\alpha ,-1/2)}}$; thus the even case of 
part a) is established.\\
Concerning the odd case, (\ref{transodd}) and (\ref{normodd}) give 
analogously
\begin{eqnarray*}
\| T_m f_o\| _{L^p_\alpha } & \approx & 
 \| \sum _{k=0}^\infty m_{2k+1} B_k h_k^{(\alpha ,1/2)} \,  
\Big( \cos \frac{\theta }{2}\Big) ^{1-2/p} \, R_k^{(\alpha ,1/2)} 
(\cos \theta ) \| _{L^p_{(\alpha ,1/2)}}\\
{} & \le & C \|  m_o \| _{M^p_{(\alpha ,1/2);(\alpha ,(p-1)/2)}}
\| f_o (\theta /2) (\cos (\theta /2))^{-2/p}\| _{L^p_{(\alpha ,1/2)}} \\
{} & \approx & 
\|  m_o \| _{M^p_{(\alpha ,1/2);(\alpha ,(p-1)/2)}}
\| f_o \| _{L^p_\alpha }, 
\end{eqnarray*}
thus
$$ \|  m \|_{M^p_{\alpha ;\alpha } \big| _{{\rm odd}}} \le C 
\|  m_o \|_{M^p_{(\alpha ,1/2);(\alpha ,(p-1)/2)} } \, .$$
The converse inequality is shown along the same lines, thus Theorem 1.1 a) is 
established.

\medskip \noindent
Concerning part b) we apply Muckenhoupt's transplantation theorem 
\cite[p. 4]{mu} twice to obtain for any sequence $\{ c_k\} $ of compact 
support, $\alpha \ge -1/2,$ and $1<p<\infty $, that 
$$ \| \sum _{k=0}^\infty c_k h_k^{(\alpha ,-1/2)} R_k^{(\alpha ,-1/2)} (\cos 
\theta )\| _{L^p_{(\alpha , -1/2)}} \approx 
\| \sum _{k=0}^\infty c_k h_k^{(\alpha ,1/2)} R_k^{(\alpha ,1/2)} (\cos \theta 
)\| _{L^p_{(\alpha ,(p -1)/2)}} $$
which in particular implies 
$$ M^p_{(\alpha ,-1/2);(\alpha ,-1/2)} =
M^p_{(\alpha ,1/2);(\alpha ,(p-1)/2)} \, ,\quad 1<p<\infty  .$$
A combination of (\ref{decompmult}) with part a) now gives part b).

\bigskip \noindent
\section{Necessary conditions for Jacobi multipliers}
Here we give a second proof of Corollary 1.2 b) which has the advantage that it 
also gives an extension to the general Jacobi case $(\alpha ,\beta ),\; 
-1/2 < \beta \le \alpha $; we note that the case $\beta = -1/2$ has already 
been discussed in \cite{hausdorff}.  
On account of the duality $M^p_{(\alpha ,\beta );(\alpha ,\beta )}=
M^{p'}_{(\alpha ,\beta );(\alpha ,\beta )},\; 1<p<\infty ,$ we can restrict 
ourselves to the case $1<p<2$ without loss of generality in the following 
(the case $p=2$ is trivial).

\thm{
Let $-1/2 < \beta \le \alpha \, ,\; 1<p<2, \; \nu =(2\beta +1)(1/p-1/2),$  and 
$\mu + \nu = (2\alpha +1)(1/p-1/2)$.
\begin{itemize}
\item[a)] If $f\in L^p_{(\alpha ,\beta )}$ is a cosine polynomial, then for some 
constant $C$ independent of $f$ there holds 
$$  \Big( \sum _{k=0}^\infty |\Delta _2^\nu \Delta ^\mu 
( s_k^{(\alpha ,\beta )}{\hat f}_{(\alpha ,\beta)} (k))|^{p'} \Big) ^{1/p'} 
\le C \| \sum _{k=0}^\infty  h_k^{(\alpha ,\beta )}{\hat f}_{(\alpha ,\beta)} 
(k) R_k^{(\alpha ,\beta )} (\cos \theta )
\| _{L^p_{(\alpha ,\beta )}} \; ,$$
where $s_k^{(\alpha ,\beta )} = (h_k^{(\alpha ,\beta )})^{1/2}$.
\item[b)] If $m\in M^p_{(\alpha ,\beta );(\alpha ,\beta )} \, $, then 
$$ \| m\|_\infty +\sup _{N\in {\bf N}}\Big( \sum _{k=N}^{2N}
|(k+1)^{\mu +\nu }\Delta _2^\nu \Delta ^\mu m_k |^{p'} \frac{1}{k+1} \Big) 
^{1/p'} \le C \| m \| _{M^p_{(\alpha ,\beta );(\alpha ,\beta )}} \, .$$
\end{itemize}
}

\bigskip \noindent
{\bf Remark} 4. Part b) is a nearly best possible necessary multiplier 
condition: one regains (up to the critical index) the right 
unboundedness domain for the Ces\`aro means (see also the following remark); 
but this example leaves open the possiblity to increase $\nu $ at the expense 
of the exponent $p',$ which would have to be replaced by some $q<p'$ 
(Sobolev embedding). That this is not possible is shown by the further example 
$m=\{ i^k (k+1)^{-\sigma } \} $. An application of part b) to $m$ 
(with $p'$ replaced by $q\le p'$) yields that 
$m$ cannot generate a bounded operator on 
$L^p_\alpha $ if $p< (2\alpha +1)/(\sigma + \alpha + 1/2),$ which 
coincides with a result of Askey and Wainger \cite[Theorem 4, ii)]{aswa}. 
There it is also proved that $m$ generates a bounded operator when 
$(2\alpha +1)/(\sigma + \alpha +1/2)<p \le 2$.

\medskip \noindent
5.  We recall that the particular case 
$\alpha =a>-1/2,\; \beta =b=-1/2,\; p < (4\alpha +4)/(2\alpha +3) $ of 
the general Cohen type inequality for Jacobi 
multipliers due to Dreseler and Soardi \cite{drso} is an immediate consequence 
of formula (3.8) in \cite{hausdorff}. So it is not surprising that 
 Theorem 3.1 b) in the case of a finite sequence $\{ m_k \} _{k=0}^N$ now 
implies the corresponding result for $-1/2 < \beta \le \alpha $. Obviously, the 
dyadic sum in part b) can be estimated from below by the single term $k=N$ and  
$\Delta _2^\nu \Delta ^\mu m_N =m_N$. A computation of the occurring $(N+1)$ 
powers immediately leads to 
\begin{equation}\label{cohen}
 (N+1)^{(2\alpha +2)(1/p-1/2)-1/2}|m_N| \le C 
\| m\| _{M^p_{(\alpha ,\beta );(\alpha ,\beta )}}\, , \quad 1<p < (4\alpha 
+4)/(2\alpha +3).
\end{equation}
We mention that, by a different method, Kalne\u\i{} \cite{kal} has obtained 
a lower bound for finite sequences in the case $p=1,\; \alpha >-1/2,\;  
\alpha \ge \beta >- 1,$ which even reflects logarithmic divergence and in 
particular implies the missing case $p=1$ in (\ref{cohen}). Kalne\u\i{}'s
lower bound is of different type than the one given in Theorem 3.1 b).

\bigskip \noindent
{\bf Proof.} First we note (cf. \cite[p. 249]{hausdorff}) that for $\mu \ge 0$
and $0 \le \theta \le \pi$
$$\Delta ^\mu \cos k\theta =\frac{1}{2} \sum _{j=0}^\infty A_j^{-\mu -1} \big(
e^{i(k+j)\theta }+ e^{-i(k+j)\theta } \big) $$
$$ = \frac12 \big( e^{ik\theta } (1-e^{i\theta })^\mu +
e^{-ik\theta } (1-e^{-i\theta })^\mu \big)
=( 2 \sin \frac{\theta }{2} ) ^\mu \cos \big( (k+\mu /2) \theta -\mu \pi /2\big) $$
$$=( 2 \sin \frac{\theta }{2} ) ^\mu \{ \cos  (k+\mu /2) \theta \cos \mu \pi /2 +
\sin (k+\mu /2) \theta \sin \mu \pi /2\} .$$
Analogously it follows that for $\nu \ge 0$ there holds
$$\Delta _2^\nu \cos (k+\mu /2)\theta 
=( 2 \sin \theta ) ^\nu \cos \big( (k+\mu /2+\nu ) \theta -\nu \pi /2\big) ,$$
$$\Delta _2^\nu \sin (k+\mu /2)\theta 
=( 2 \sin \theta ) ^\nu \sin \big( (k+\mu /2+\nu ) \theta -\nu \pi /2\big) .$$
Hence we obtain the following $(L^1,l^\infty )$--estimate for a trigonometric 
polynomial $f$
\begin{equation}\label{L1estimate}
|\Delta _2^\nu \Delta ^\mu a_k|
\le C \int _0^\pi |f(\theta )| \Big( \sin \frac{\theta }{2} \Big)^{\mu + \nu }
\Big( \cos \frac{\theta }{2}\Big) ^\nu \, d\theta \, ,
\end{equation} 
where $a_k=\int _0^\pi f(\theta ) \cos k\theta \, d\theta .$
For a corresponding $(L^2,l^2)$--estimate we observe that 
$$\Delta _2^\nu \Delta ^\mu a_k 
= \int _0^\pi f(\theta ) \Big( \sin \frac{\theta }{2} \Big)^{\mu + \nu }
\Big( \cos \frac{\theta }{2}\Big) ^\nu \big( g_1(\theta )\cos k\theta +g_2(\theta ) 
\sin k\theta \big)\, d\theta \, =: C_k +D_k,$$
where the continuous functions $g_l$ are linear combinations of $\sin c\theta $ 
and $\cos c \theta ,\; c$ denoting different constants depending only upon $\mu $ 
and $\nu $. If we now consider functions $f$ with 
$f(\theta ) \Big( \sin \frac{\theta }{2} 
\Big)^{\mu + \nu }\Big( \cos \frac{\theta }{2}\Big) ^\nu  \in L^2(0,\pi )$ and 
observe that the systems $\{ \sin k\theta \} $ and $\{ \cos k\theta \} $ are 
essentially orthonormal, it follows by the Parseval formula that
\begin{equation}\label{L2estimate}
\sum _{k=0}^\infty |\Delta _2^\nu \Delta ^\mu a_k|^2  \le C \left(
\sum _{k=0}^\infty |C_k|^2 + \sum _{k=0}^\infty |D_k|^2 \right) \le C 
\int _0^\pi |f(\theta ) \Big( \sin \frac{\theta }{2} 
\Big)^{\mu + \nu }\Big( \cos \frac{\theta }{2}\Big) ^\nu  |^2d\theta \, .
\end{equation}
The Riesz--Thorin interpolation theorem applied to (\ref{L1estimate}) and 
(\ref{L2estimate}) gives for $1\le p \le 2$ the Hausdorff--Young type inequality
\begin{equation}\label{haus}
\Big( \sum _{k=0}^\infty |\Delta _2^\nu \Delta ^\mu a_k|^{p'} \Big) ^{1/p'}  
\le C \Big( \int _0^\pi |\sum _{k=0}^\infty a_k \cos k\theta 
\Big( \sin \frac{\theta }{2} 
\Big)^{\mu + \nu }\Big( \cos \frac{\theta }{2}\Big) ^\nu  |^p d\theta \Big) ^{1/p}\, .
\end{equation}
If one transplants this  inequality for the cosine expansion (which 
corresponds in the Jacobi setting to the parameters $(-1/2,-1/2)$ )
to arbitrary Jacobi expansions with parameters $(\alpha ,\beta ), 
-1/2 < \beta \le \alpha ,$ then in $L^p$--spaces with natural weights one has to 
check the hypotheses in Muckenhoupt's transplantation theorem \cite[p. 4]{mu}.
For $1<p\le 2$ this leads to the restrictions $\nu = (2\beta +1)(1/p-1/2)>0$ 
(hence $\beta =-1/2$ is not admitted) and $\mu +\nu >0,\; \mu = 
2(\alpha -\beta )(1/p-1/2)$. Now choose 
$$ a_k=\int _0^\pi f(\theta ) \Big( \sin \frac{\theta }{2}\Big) ^{\alpha +1/2} 
\Big( \cos \frac{\theta}{2}\Big) ^{\beta +1/2} \phi _k^{(\alpha ,\beta )}
(\theta )\, d\theta = s_k^{(\alpha ,\beta )} {\hat f} _{(\alpha ,\beta )}(k)
\,   ,$$
where we use the Muckenhoupt notation \cite[(2.2)]{mu}
$$ \phi _k^{(\alpha ,\beta )}(\theta ) = t_k^{(\alpha ,\beta )}
P_k^{(\alpha ,\beta )}(\cos \theta )\Big( \sin \frac{\theta }{2}\Big) 
^{\alpha +1/2} \Big( \cos \frac{\theta }{2}\Big) ^{\beta +1/2}$$
with $t_k^{(\alpha ,\beta )} =s_k^{(\alpha ,\beta )} /P_k^{(\alpha ,\beta )} 
(1) .$ Then Muckenhoupt's transplantation theorem gives
$$\Big( \int _0^\pi \Big| \sum _{k=0}^\infty a_k \cos k\theta \Big| ^p
\Big( \sin \frac{\theta }{2} \Big)^{p(\mu + \nu )}
\Big( \cos \frac{\theta }{2} \Big) ^{p\nu }   d\theta \Big) ^{1/p}$$
$$ \le C \Big( \int _0^\pi \Big| \sum _{k=0}^\infty a_k \phi _k^{(\alpha ,\beta 
)}(\theta )\Big| ^p \Big( \sin \frac{\theta }{2} \Big)^{p(\mu + \nu ) }
\Big( \cos \frac{\theta }{2} \Big) ^{p\nu } d\theta \Big) ^{1/p} $$
$$ =  C \Big( \int _0^\pi \Big| \sum _{k=0}^\infty {\hat f} _{(\alpha ,\beta )}
(k)  h_k^{(\alpha ,\beta )} R_k^{(\alpha ,\beta )}(\cos 
\theta )\Big| ^p \Big( \sin \frac{\theta }{2} \Big) ^{2\alpha +1}
\Big( \cos \frac{\theta }{2}\Big) ^{2\beta +1} d\theta \Big) ^{1/p}. $$
A combination with (\ref{haus}) gives part a) of Theorem 3.1.

\bigskip \noindent
Concerning part b), consider a $C^\infty $--function $\chi (x)$ with 
$$ \chi (x) = \left\{ \begin{array}{r@{\quad{\rm  if}\quad}l}
0 & 0\le x \le 1/2 \\
1 & 1\le x \le 4 \\
0 & x\ge 8
\end{array} \right. \; , \quad \chi _i(x) = \chi (2^{-i}x) \, , $$
and an associated test sequence $\{ (s_k^{(\alpha ,\beta )})^{-1} \chi _i(k)\} 
$. Then, by \cite[Theorem 2]{batr}, it is not hard to see that
\begin{equation}\label{testfctn}
 \| \sum _{k=0}^\infty (s_k^{(\alpha ,\beta )})^{-1} \chi _i(k) 
    h_k^{(\alpha ,\beta )} R_k^{(\alpha ,\beta )}(\cos \theta ) \| _{L^p
_{(\alpha ,\beta )}}  \le C (2^i)^{(2\alpha +2)/p' -\alpha -1/2}\, .
\end{equation}
By part a) and the hypothesis $m\in  M^p_{(\alpha ,\beta );(\alpha ,\beta )}$ 
we have that
$$ \Big( \sum _{k=2^i}^{2^{i+1}}
|\Delta _2^\nu \Delta ^\mu (m_k \chi _i(k))|^{p'} \Big) ^{1/p'}  \le  
\Big( \sum _{k=0}^\infty |\Delta _2^\nu \Delta ^\mu (m_k \chi _i(k)) |^{p'}
\Big) ^{1/p'} \quad \quad \quad \quad \quad \quad \quad \quad $$
$$ \quad \quad \quad \quad 
 \le C  \| \sum _{k=0}^\infty  m_k (s_k^{(\alpha ,\beta )})^{-1} 
\chi _i(k) h_k^{(\alpha ,\beta )} R_k^{(\alpha ,\beta )} (\cos \theta )
\| _{L^p_{(\alpha ,\beta )}} $$
$$ \quad \quad \quad \quad \quad \quad 
 \le C  \| m\| _{M^p_{(\alpha ,\beta );(\alpha ,\beta )}}
\| \sum _{k=0}^\infty   (s_k^{(\alpha ,\beta )})^{-1} \chi _i(k) 
h_k^{(\alpha ,\beta )} R_k^{(\alpha ,\beta )} (\cos \theta )
\| _{L^p_{(\alpha ,\beta )}} $$
whence by (\ref{testfctn})
$$\Big( \sum _{k=2^i}^{2^{i+1}}
|(k+1)^{\mu +\nu }\Delta _2^\nu \Delta ^\mu (m_k \chi _i(k)) |^{p'} 
\frac{1}{k+1} \Big) ^{1/p'} 
\le C \| m \| _{M^p_{(\alpha ,\beta );(\alpha ,\beta )}} $$
with the right side independent of $i$. The final statement b) now follows 
along the lines of the proof of \cite[Lemma 2.3]{gtlag}.

\bigskip \noindent
{\bf Remark} 6. Of course one can state analogous to part b) the same 
necessary conditions for the above considered cosine expansions in 
weighted $L^p$--spaces by applying once more Muckenhoupt's transplantation 
theorem. Observe that then multipliers only make sense for those functions 
whose first $N$ coefficients of the cosine expansion vanish; here 
$N= {\rm max}\{ [(\alpha +1/2)(1/p-1/2)+1/2p +1/2],0\} +
{\rm max}\{ [(\beta +1/2)(1/p-1/2)+1/2p +1/2],0\} $ .

\bigskip \noindent
\section{Criteria for integrable functions}
First we consider the problem:  Given a sequence $\{ f_k\} $, what are 
sufficient conditions satisfied by $\{ f_k\} $ such that the $f_k$ are 
Fourier--Jacobi coefficients of an $L^1_\alpha $--function $f$?
Via the transformation formulas in (\ref{quadratic}) we briefly 
give improvements of known criteria. We start with a Parseval relation.
\prop{
Fix $\alpha \ge -1/2$ and let $f(\theta )=\sum {\hat f}_{(\alpha ,\alpha )}(k) h_k 
^{(\alpha ,\alpha )} R_k ^{(\alpha ,\alpha )} (\cos \theta )$ be a finite sum 
(i.e., a polynomial in $\cos \theta $).
\begin{itemize}
\item[a)] If $-1/2 < \mu < \alpha +2, $ then
$$ \int _0^\pi |f(\theta )|^2 (\sin \theta )^{2(\alpha +\mu )+1}d\theta 
\le C \sum _{k=0}^\infty |\Delta _2^\mu {\hat f}_{(\alpha ,\alpha )}(k)|^2 h_k 
^{(\alpha ,\alpha )} \, .$$
\item[b)] If $\mu > -1,$ then the converse holds, i.e.
$$ \sum _{k=0}^\infty |\Delta _2^\mu {\hat f}_{(\alpha ,\alpha )}(k)|^2 h_k 
^{(\alpha ,\alpha )} \le C 
\int _0^\pi |f(\theta )|^2 (\sin \theta )^{2(\alpha +\mu )+1}d\theta \, .$$
\end{itemize}
}
For the proof we have only to observe that by (\ref{transeven}) and 
(\ref{transodd}) 
$$\int _0^\pi |f(\theta )|^2 (\sin \theta )^{2(\alpha +\mu )+1} d\theta 
\approx \int _0^{\pi /2}\Big( |f_e(\theta )|^2 +|f_o(\theta )|^2 \Big) 
(\sin \theta )^{2(\alpha +\mu )+1} d\theta $$
$$\approx \int _0^\pi \Big| \sum _{k=0}^\infty A_k h_k^{(\alpha ,-1/2)} R_k^{(
\alpha ,-1/2)}(\cos \theta )\Big| ^2 \Big( \sin \frac{\theta }{2}\Big) 
^{2(\alpha +\mu )+1} d\theta $$
$$+\int _0^\pi \Big| \sum _{k=0}^\infty B_k h_k^{(\alpha ,1/2)} 
R_k^{(\alpha ,1/2)}(\cos \theta )\Big| ^2 \Big( \sin \frac{\theta }{2}\Big) 
^{2(\alpha +\mu )+1} \Big( \cos \frac{\theta }{2}\Big) ^2 d\theta \, .$$
Then \cite[Theorem 1]{williamstown}, whose proof extends to the case $\alpha 
,\beta \ge -1/2$, can be applied to the two terms of the 
right side, and the assertion follows after noting that
$$\sum _{k=0}^\infty |\Delta ^\mu A_k|^2 h_k^{(\alpha ,-1/2)} +
\sum _{k=0}^\infty |\Delta ^\mu B_k|^2 h_k^{(\alpha ,1/2)} 
\approx \sum _{k=0}^\infty |\Delta _2^\mu {\hat f}_{(\alpha ,\alpha )}(k)|^2 
h_k^{(\alpha ,\alpha )} \, .$$

\bigskip \noindent
\thm{
Let $\alpha \ge -1/2$ and $\mu >\alpha +1$.
If $\{ c_k\}$ is a bounded sequence with $\lim _{k\to \infty } c_k=0$ and 
$$ \sum _{j=1}^\infty \Big( \sum _{k=2^{j-1}}^{2^j-1} k^{-1}|c_k|^2 \Big) 
^{1/2} + \sum _{j=1}^\infty 
\Big( \sum _{k=2^{j-1}}^{2^j-1} k^{-1}|k^\mu \Delta _2^\mu c_k |^2 
\Big) ^{1/2} \le K_{ \{ c_k\} } \, ,$$
then there exists an $f\in L^1_\alpha $ with 
${\hat f}_{(\alpha ,\alpha )}(k)=c_k$ for all $k\in {\bf N}_0$ and $\| f\| 
_{L^1_\alpha } \le C  K_{ \{ c_k\} } \, .$
}

\medskip \noindent
The proof follows from Proposition 4.1 a), analogous to that of 
\cite[Theorem 2 a]{williamstown}, or directly from 
\cite[Theorem 2 a]{williamstown} by the same method used for Proposition 4.1.

\medskip \noindent
Next we give another simple sufficient multiplier condition which is not 
comparable with Theorem 4.2 --- see the discussion in \cite{williamstown}.

\bigskip \noindent
\thm{
Let $\alpha \ge -1/2$ and $\mu >\alpha +1/2$.
If $\{ c_k\}$ is a bounded sequence with $\lim _{k\to \infty } c_k=0$ and 
$$ \sum _{k=0}^\infty (k+1)^\mu |\Delta _2^{\mu +1} c_k| 
\le K_{ \{ c_k\} } \, ,$$
then there exists an $f\in L^1_\alpha $ with 
${\hat f}_{(\alpha ,\alpha )}(k)=c_k$ for all $k\in {\bf N}_0$ and $\| f\| 
_{L^1_\alpha } \le C  K_{ \{ c_k\} } \, .$
}

\medskip \noindent
Split the sequence $\{ c_k\} $ into the two subsequences $\{ a_k\} $
and $\{ b_k\} $, where $a_k= c_{2k}\, ,\; b_k= c_{2k+1} $, and 
observe that
$$\sum _{k=0}^\infty  (k+1)^\mu |\Delta _2^{\mu +1} c_k| \approx 
\sum _{k=0}^\infty  (k+1)^\mu |\Delta ^{\mu +1} a_k| +
\sum _{k=0}^\infty  (k+1)^\mu |\Delta ^{\mu +1} b_k|.$$
Now one can follow for each 
subsequence the proof of \cite[Lemma 1]{hausdorff}.  
The assumption there that the sequence has compact 
support is not used in \cite[(3.5)]{hausdorff}. First observe that 
$$ \int _0^\pi \Big| \sum _{j=0}^k (A^\mu _{k-j} /A^\mu _k) 
h_j^{(\alpha ,-1/2)} R_j^{(\alpha ,-1/2)}(\cos \theta )
\Big| \Big( \sin \frac{\theta }{2} \Big) 
^{2\alpha +1} d\theta \le K,\quad \mu >\alpha +1/2, $$
which is proved in \cite[(9.41.1)]{szego}. Thus the series 
$$  \sum _{k=0}^\infty  A_k^\mu \Delta ^{\mu +1} a_k 
\sum _{j=0}^k (A^\mu _{k-j} /A^\mu _k) h_j^{(\alpha ,-1/2)}
R_j^{(\alpha ,-1/2)}(\cos \theta ) $$
converges a.e. to a function $f_1 \in L^1_{(\alpha ,-1/2)}$ with 
coefficients $(f_1)\, {\hat {}}_{(\alpha ,-1/2)}(k) = a_k=c_{2k}$. Analogously 
one deals with the sequence $\{ b_k\} $ to which one associates the function
$$ f_2(\theta ) =  \sum _{k=0}^\infty  A_k^\mu \Delta ^{\mu +1} b_k 
\sum _{j=0}^k (A^\mu _{k-j} /A^\mu _k) h_j^{(\alpha ,1/2)}
R_j^{(\alpha ,1/2)}(\cos \theta ). $$
To deduce that $f_2\in L^1_{(\alpha , 0)}$ 
one needs the following boundedness result concerning the Ces\`aro kernel 
$$ \int _0^\pi \Big| \sum _{j=0}^k (A^\mu _{k-j} /A^\mu _k) h_j^{(\alpha ,1/2)}
R_j^{(\alpha ,1/2)}(\cos \theta )\Big| \Big( \sin \frac{\theta }{2} \Big) ^{2
\alpha +1} \cos \frac{\theta}{2} \, \, d\theta \le K,\quad \mu >\alpha +1/2, $$
which follows by a slight modification of Szeg\"o's proof --- note that by  
the third case of \cite[(7.34.1)]{szego},  
$\, \int _{\pi /2}^\pi |P_n^{(\alpha ,1/2)}(\cos \theta )| \cos \theta /2 \, 
d\theta = O(n^{-1/2})$, so that the right side estimate in 
\cite[(9.41.2)]{szego}) remains valid, as does the rest of the proof in 
\cite{szego}.  To complete the proof it suffices to set $f=f_e+f_o$ with 
$f_e(\theta ) =f_1(2\theta )$ and $f_o(\theta )= \cos \theta \, f_2(2\theta ),$ 
and to use (\ref{decompnorm}) and (\ref{normeven}) -- (\ref{transodd}).

\bigskip \noindent
Let us turn to the question of necessary conditions. As in Sec.~1 
decompose a co\-sine polynomial $f\in L^1_\alpha $ into its even and odd parts
with respect to the line $\theta =\pi /2$: 
$f=f_e+f_o$, and set $f_1(\theta )=f_e(\theta /2)$. 
Then, by (\ref{transcoeffeven}), $A_k:= (f_1)\, {\hat {}}_{(\alpha ,-1/2)}(k) 
= 2^{2\alpha+1} (f_e)\, {\hat {}}_{(\alpha ,
\alpha )}(2k),$ and \cite[(3.2)]{hausdorff} gives for $\alpha \ge -1/2,\; 
\nu \ge 0$,  that 
$$ \sup_k |(h_{k}^{(\alpha ,-1/2)})^{1/2} \Delta^\nu A_k| \approx 
\sup_k |(h_{2k}^{(\alpha ,\alpha )})^{1/2} \Delta_2^\nu
(f_e)\, {\hat {}}_{(\alpha ,\alpha)}(2k)| $$
$$\quad \quad \quad \quad \quad \quad \quad \quad \le C 
 \int _0^\pi | (\sin \theta )^{\alpha +\nu +1/2}f_e(\theta )| d\theta.$$
Similarly, from the case $\alpha \ge -1/2, \beta = 1/2$ of
\cite[(3.2)]{hausdorff} (it extends immediately to this case) 
applied to the function
$f_2(\theta) = f_o(\theta/2)/\cos (\theta/2)$ and the fact that
$B_k:=(f_2)\, {\hat {}}_{(\alpha ,1/2)}(k) =  2^{2\alpha+1}
(f_o)\, {\hat {}}_{(\alpha , \alpha)}(2k+1)$
we get
$$ \sup_k |(h_k^{(\alpha ,1/2)})^{1/2} \Delta ^\nu B_k | \approx 
\sup_k |(h_{2k+1}^{(\alpha ,\alpha )})^{1/2} \Delta _2^\nu 
(f_o)\, {\hat {}}_{(\alpha ,\alpha)}(2k+1)| $$
$$ \quad \quad \quad \quad \quad \quad \quad \quad \le C 
 \int _0^\pi | (\sin \theta )^{\alpha +\nu+1/2}f_o(\theta )|\,  d\theta .$$
\noindent
Combining these two estimates gives the inequality
\begin{equation}\label{supinequality}
\sup_k |(h_k^{(\alpha ,\alpha )})^{1/2} \Delta_2^\nu
{\hat f}_{(\alpha ,\alpha )}(k)| \le C 
 \int _0^\pi | (\sin \theta )^{\alpha +\nu+1/2}f(\theta )| d\theta.
\end{equation}
\noindent
Since $h_k^{(\alpha ,\alpha )} \approx (k+1)^{2\alpha +1,},$ application
of the Riesz--Thorin theorem to (\ref{supinequality}) and 
Proposition 4.1 b) yields part a) of

\smallskip \noindent
\thm{ 
a) Let $1 \le p \le 2, \alpha \ge -1/2,$ and $ \nu \ge 0.$ 
If $f$ is a cosine polynomial, then the Hausdorff--Young type
inequality
$$\Big( \sum_{k=0}^\infty  |(k+1)^{\alpha +1/2} \Delta_2^{\nu}
{\hat f}_{(\alpha ,\alpha )}(k)|^{p'} \Big)^{1/p'} 
\le C \Big( \int _0^\pi | (\sin \theta )^{\alpha +\nu+1/2}
f(\theta )|^p\,  d\theta \Big)^{1/p}$$
holds.  Also, if  $f\in L^1_\alpha $ then
$$\sup _{k\in {\bf N}_0} |(k+1)^{\alpha +1/2} \Delta_2^{\alpha +1/2}
{\hat f}_{(\alpha ,\alpha )}(k)| \le C \| f\|_{L^1_\alpha } \, ,$$
which gives a necessary condition for a sequence to be the sequence
of Fourier--Jacobi coefficients of an $L^1_\alpha $--function. 

\medskip \noindent
b) If $0<\nu < \alpha +1/2, $ then
$$ \sum _{k=0}^\infty (k+1)^{\nu -1}|\Delta _2 ^\nu  
{\hat f}_{(\alpha ,\alpha )}(k)| \le C \| f\|_{L^1_\alpha } .$$
}

\medskip \noindent
{\bf Remark} 7. The case $p=1$ of Part a) contains a Cohen type inequality for 
ultraspherical expansions. The assertion in part b) does not follow from part 
a): Observe that in general for $0<\nu < \alpha +1/2$ there only holds
$$\sup _{k\in {\bf N}_0} |(k+1)^\nu \Delta_2^\nu 
{\hat f}_{(\alpha ,\alpha )}(k)| \le C 
\sup _{k\in {\bf N}_0} |(k+1)^{\alpha +1/2} \Delta_2^{\alpha +1/2}
{\hat f}_{(\alpha ,\alpha )}(k)| $$
(consider e.g. the sequence $\{ k^{i\gamma } \} \, ,\gamma \in {\bf R},$ 
fixed); then it is clear that the estimate of part a) would lead to the 
diverging harmonic series $\sum 1/(k+1)$.

\medskip \noindent
The proof of Theorem 4.4 b) is an immediate consequence of 
\la{
Let $\alpha > -1/2,\; 0< \nu <\alpha +1/2$, and $ 0 < \theta < \pi /2$. 
Then there holds\\
$a) \quad \quad \quad \quad \quad \quad \quad |\Delta _2^\nu R_k^{(\alpha ,
\alpha )}(\cos \theta )| \le C \, (\sin \theta )^\nu \, ,$

\smallskip \noindent
$b) \quad \quad \quad \quad |\Delta _2^\nu R_k^{(\alpha ,\alpha )}(\cos \theta )
| \le C \, (\sin \theta )^{\nu -\alpha -1/2} (k+1)^{-\alpha -1/2},$

\smallskip \noindent
where in b) it is additionally assumed that $\nu \ge [\alpha +1/2]^* $ if 
$\alpha >1/2$. Here we use the notation 
$[a]^*,\; a\in {\bf R},$ for the greatest integer smaller than 
$a,\; [a]^*<a.$
}

\bigskip \noindent
For suppose that Lemma 4.5 is true. Obviously,
$$ \sum _{k=0}^\infty (k+1)^{\nu -1}|\Delta _2 ^\nu  
{\hat f}_{(\alpha ,\alpha )}(k)| \le C \| f\|_{L^1_\alpha } \sup _{0\le \theta 
\le \pi } \sum _{k=0}^\infty (k+1)^{\nu -1}|\Delta _2^\nu  
R_k^{(\alpha ,\alpha )}(\cos \theta )| \, .$$
On account of the symmetry of the ultraspherical polynomials (in each 
difference all the polynomials are even or all are odd with respect to the 
line $\theta =\pi /2$) we may take the supremum over $0\le \theta \le \pi /2$.
Now decompose the interval $[0,\pi /2]$  into intervals $I_j:=[2^{-j-1}\pi ,
2^{-j}\pi ]\, ,\; j\in {\bf N},$ and consider $\theta \in I_j$. If we set 
$N_j=[2^j/ \pi ]$, then $0< C' \le N_j \sin \theta \le C < \infty $ uniformly 
in $j$. It then follows by the above two estimates that
$$\sup _{\theta \in I_j} \sum _{k=0}^\infty (k+1)^{\nu -1} |\Delta _2^\nu  
R_k^{(\alpha ,\alpha )}(\cos \theta )| \quad \quad \quad \quad \quad \quad \quad 
\quad \quad \quad \quad \quad $$
$$\le C ( \sin \theta )^\nu \sum _{k=0}^{N_j} (k+1)^{\nu -1} + C 
(\sin \theta )^{\nu -\alpha -1/2} \sum _{k=N_j}^\infty (k+1)^{\nu -\alpha -3/2}
\le C .$$
The convergence of the last series is ensured by the hypothesis $\nu < \alpha 
+1/2 $; thus the assertion holds for $\alpha _0:=[\alpha +1/2]^* \le \nu < 
\alpha +1/2$. The extension to all $\nu, \; 0< \nu <\alpha +1/2, $ is
straightforward. It is clear that 
$|{\hat f}_{(\alpha ,\alpha )}(k)| \le  \| f\|_{L^1_\alpha } $ since 
$|R_k^{(\alpha ,\alpha )}(\cos \theta )| \le 1 $ for $\alpha >-1/2$. 
Thus proceeding as in \cite[Lemma 1]{wbv}, in particular using the Andersen  
formula for bounded sequences $\{ a_k\} $
$$\Delta ^{\lambda +\kappa } a_k = \Delta ^\lambda (\Delta ^\kappa  a_k), \quad 
\kappa \ge 0,\; \lambda >-1,\; \kappa +\lambda >0,$$
one obtains for ${\rm min} \{ 0,\, \alpha _0 -1\} < \nu <\alpha _0 $
$$ \sum_{k=0}^\infty (k+1)^{\nu -1} |\Delta ^\nu  {\hat f}_{(\alpha ,\alpha )}
(k)| \le  C \sum_{k=0}^\infty A_k^{\nu -1} |\sum _{j=k}^\infty A_{j-k}^{\alpha 
_0 -\nu -1} \Delta ^{\alpha _0} {\hat f}_{(\alpha ,\alpha )}(j)|$$
$$ \le C \sum_{j=0}^\infty |\Delta ^{\alpha _0} {\hat f}_{(\alpha ,\alpha )}(j)|
\sum _{k=0}^j A_k^{\nu -1} A_{j-k}^{\alpha _0 -\nu -1} \le C \sum_{j=0}^\infty 
(k+1)^{\alpha _0 -1}|\Delta ^{\alpha _0} {\hat f}_{(\alpha ,\alpha )}(j)|\, .$$
Iteration of this procedure gives the assertion.

\bigskip \noindent
{\bf Proof of Lemma 4.5.} Mehler's integral \cite[10.9 (32)]{htf}  and the 
formula for the fractional difference $\Delta _2^\nu \cos (k+\mu /2) 
\theta $ preceding (\ref{L1estimate}) give
$$ \Delta _2^\nu R_k^{(\alpha ,\alpha )}(\cos \theta ) \quad \quad \quad \quad 
\quad \quad \quad \quad \quad \quad \quad \quad \quad \quad \quad \quad \quad 
\quad \quad \quad \quad \quad $$
$$=C_\alpha (\sin \theta )^{-2\alpha} \int _0^\theta  (\cos \phi -\cos \theta 
)^{\alpha -1/2} (2\sin \phi )^\nu \cos ((k+\alpha +\nu +1/2)\phi -\nu \pi /2 
)\, d\phi \, .$$
Hence 
$$ |\Delta _2^\nu R_k^{(\alpha ,\alpha )}(\cos \theta ) | 
 \le C (\sin \theta )^{-2\alpha} \int _0^\theta  \Big( \sin \frac{\theta 
+\phi }{2} \, \sin \frac{\theta -\phi }{2} \Big) ^{\alpha -1/2} (\sin \phi 
)^\nu d\phi $$
$$ \le C (\sin \theta )^{-\alpha -1/2} \bigg[ \int _0^{\theta /2} + 
\int _{\theta /2} ^\theta \bigg] \Big( \sin \frac{\theta -\phi }{2} \Big) 
^{\alpha -1/2} (\sin \phi )^\nu d\phi $$
$$\le C (\sin \theta )^{-1} \int _0^{\theta /2} \phi ^\nu d\phi 
+C (\sin \theta )^{\nu -\alpha -1/2} \int _{\theta /2} ^\theta (\theta -\phi 
)^{\alpha -1/2} d\phi \le C (\sin \theta )^\nu $$
since $\alpha >-1/2$ and $0< \theta < \pi /2$; thus part a) is established.

\medskip \noindent
The case $\nu \in {\bf N}$ of part b) has already been shown by 
Kalne\u\i{} \cite[Lemma 3]{kal}. Since it is clear by part a) that b) holds for 
$0<\theta \le c/(k+1)$ for fixed $c>0$, without loss of generality we 
may assume $\pi /k<\theta \le \pi /2,\; k\ge 3$.  Obviously,
$$ \Delta _2^\nu R_k^{(\alpha ,\alpha )}(\cos \theta ) =C 
(\sin \theta )^{\nu -\alpha -1/2} I_{\alpha ,\nu} (\theta ;k) ,$$
where
$$I_{\alpha ,\nu} (\theta ;k) = \int _0^\theta  \Big( \frac{
\cos \phi -\cos \theta }{\sin \theta } \Big) ^{\alpha -1/2} \Big( \frac{
\sin \phi }{\sin \theta} \Big) ^\nu \cos ((k+\alpha +\nu +1/2)\phi -\nu \pi /2 
)\, d\phi \, .$$
Now write 
$$ \cos ((k+\alpha +\nu +1/2)\phi -\nu \pi /2 ) \quad \quad \quad \quad \quad 
\quad \quad \quad \quad \quad \quad \quad \quad \quad \quad \quad \quad \quad $$
$$= \cos ((\alpha + \nu +1/2)\phi -\nu \pi /2 ) \, \cos k\phi  -
\sin ((\alpha +\nu +1/2)\phi -\nu \pi /2 ) \, \sin k\phi \, .$$
The idea for obtaining the $(k+1)^{-\alpha -1/2}$ decrease is to interpret the 
preceding integral as cosine and sine coefficients of functions which satisfy 
appropriate $L^1$--Lipschitz conditions. Then formula (4.2) in \cite[Chap. 
II]{zy} and iterated integrations by parts of sufficiently high order give the 
desired $(k+1)$--decrease. 

\smallskip \noindent
Let us first look at the case $-1/2 <\alpha \le 1/2$ and set
$$ G_{\alpha ,\nu ;\theta }(\phi ) = 
\left\{ \begin{array}{l@{\quad{\rm if}\quad}l}
\Big( \frac{
\cos \phi -\cos \theta }{\sin \theta } \Big) ^{\alpha -1/2} \Big( \frac{
\sin \phi }{\sin \theta} \Big) ^\nu  & 0 \le \phi \le \theta \le \pi /2 \\
0 & \theta < \phi <\pi \, .
\end{array} \right.  $$
Then 
\begin{eqnarray*}
I_{\alpha ,\nu} (\theta ;k) & = & \int _0^\pi G_{\alpha ,\nu ;\theta }(\phi ) 
\cos ((\alpha + \nu +1/2)\phi -\nu \pi /2 ) \cos k\phi \, d\phi \\
{} & {} & \quad - \int _0^\pi G_{\alpha ,\nu ;\theta }(\phi ) 
\sin ((\alpha + \nu +1/2)\phi -\nu \pi /2 ) \sin k\phi \, d\phi \, .
\end{eqnarray*}
Since $\cos ((\alpha + \nu +1/2)\phi -\nu \pi /2 ) $ and $\sin ((\alpha + \nu +
1/2)\phi -\nu \pi /2 ) $ are bounded $C^\infty $--functions we can clearly 
neglect them when discussing the smoothness of $G_{\alpha ,\nu ;\theta }(\phi )
  $. Elementary, though tedious, computations show that 
$$ \int _0^\pi |G_{\alpha ,\nu ;\theta }(\phi +\delta ) -
G_{\alpha ,\nu ;\theta }(\phi ) |\, d\phi \le C 
\delta ^{\alpha +1/2},\quad 0< \delta \le \theta ,\quad \theta \ge \pi /k\, .$$
Now formula (4.2) in \cite[Chap. II]{zy} (adapted for sine and cosine 
expansions) gives 
$$ I_{\alpha ,\nu} (\theta ;k) \le C (k+1)^{-\alpha - 1/2}\, ,$$ 
which yields assertion b) in the case $-1/2 < \alpha \le 1/2$. If $1/2 < \alpha \le 
3/2$, an integration by parts leads to 
$$ I_{\alpha ,\nu} (\theta ;k) = -\frac{1}{k+\alpha + \nu +1/2}
\int _0^\theta G'_{\alpha ,\nu ;\theta }(\phi ) 
\sin ((k+\alpha + \nu +1/2)\phi -\nu \pi /2 ) \, d\phi .$$
An examination of the derivative $G'_{\alpha ,\nu ;\theta }(\phi ) $ 
shows that it has at least the same smooth\-ness as $G_{\alpha -1 ,\nu ;\theta }
(\phi ) $. Hence again $ I_{\alpha ,\nu} (\theta ;k) \le C 
(k+1)^{-\alpha - 1/2}\, $. An iteration of this procedure finally shows the 
assertion b) of Lemma 4.5 to be true for all $\alpha > -1/2$.

\end{document}

%% file: macrola2.tex
\newtheorem{THEOREM}{Theorem}[section]
\newcommand{\thm}[1]{
\begin{THEOREM}
#1
\end{THEOREM}
}
\newtheorem{PROPOSITION}[THEOREM]{Proposition}
\newcommand{\prop}[1]{
\begin{PROPOSITION}
#1
\end{PROPOSITION}
}
\newtheorem{rem}{Remark}[section]

\newtheorem{LEMMA}[THEOREM]{Lemma}
\newcommand{\la}[1]{
\begin{LEMMA}
#1
\end{LEMMA}
} 
\newtheorem{COROLLARY}[THEOREM]{Corollary}
\newcommand{\coro}[1]{
\begin{COROLLARY}
#1
\end{COROLLARY}
}
\newtheorem{DEFINITION}{Definition}[section]
 